\documentclass[12pt]{article}

\usepackage{multirow}
\usepackage{color}
\usepackage{epsfig}
\usepackage{graphicx}
\usepackage{amsmath}
\usepackage{color}
\usepackage{pstricks}
\usepackage{amsfonts}
\usepackage{amssymb}
\usepackage{amsthm}

\numberwithin{equation}{section}

\theoremstyle{plain}
\newtheorem{thm}{Theorem}

\theoremstyle{definition}

\theoremstyle{remark}

\theoremstyle{proof}

\usepackage{ulem}
\usepackage{color}
\definecolor{purple}{rgb}{0.40,0.05,0.98} \definecolor{rougefonce}{rgb}{0.70,0.00,0.10}
\usepackage[colorlinks=true,linkcolor=purple,citecolor=rougefonce]{hyperref}
\usepackage{tikz} 
\usepackage[sort]{natbib} 

\usepackage{booktabs}
\usepackage{lscape}  
\usepackage{multicol}
\usepackage{longtable}

\DeclareMathOperator{\logit}{logit}

\begin{document}

\thispagestyle{empty}
\title{Estimation under cross-classified sampling with application to a childhood survey}
\date{\today}
\author{H\'el\`ene Juillard\thanks{INED, 133 boulevard Davout, 75020 Paris, France} \\
Guillaume Chauvet\thanks{CREST/ENSAI, Campus de Ker Lann, 35170 Bruz, France} \\
Anne Ruiz-Gazen\thanks{Toulouse School of Economics, 21 all\'ee de Brienne, 31000 Toulouse, France}}

\baselineskip=24pt

\maketitle

%
%
%
%
%


\begin{center}
{\large{\bf Abstract}}
\end{center}

The cross-classified sampling design consists in drawing samples from a two-dimension population, independently in each dimension. Such design is commonly used in consumer price index surveys and has been recently applied to draw a sample of babies in the French ELFE survey, by crossing a sample of maternity units and a sample of days. We propose to derive a general theory of estimation for this sampling design. We consider the Horvitz-Thompson estimator for a total, and show that the cross-classified design will usually result in a loss of efficiency as compared to the widespread two-stage design. We obtain the asymptotic distribution of the Horvitz-Thompson estimator, and several unbiased variance estimators. Facing the problem of possibly negative values, we propose simplified non-negative variance estimators and study their bias under a super-population model. The proposed estimators are compared for totals and ratios on simulated data. An application on real data from the ELFE survey is also presented, and we make some recommendations. Supplementary materials are available online.


\vspace{1mm}

\noindent
{\bf Some key words:} analysis of variance, Horvitz-Thompson estimator, independence, invariance, Sen-Yates-Grundy estimator, two-phase sampling, two-stage sampling.

\vspace{1mm}

\noindent
{\bf Short title}: Estimation under cross-classified sampling


\section{\baselineskip=20pt Introduction} \label{sec:Intro}
\baselineskip=24pt


\noindent The 2011 French Longitudinal Survey on Childhood (ELFE) comprises more than 18,000 children selected on the basis of their place and date of birth. On the one hand, a sample of 320 maternity units has been drawn on the metropolitan territory. On the other hand, a sample of 25 days divided in four time periods and spread almost equally across the four seasons of 2011 has been selected. The babies born at the sampled locations and on the sampled days have been approached through midwives. Data were collected on babies whose parents consented to their inclusion during their stay at the maternity. ELFE is conducted by the National Institute for Demographic Studies (INED), the National Institute for Health and Medical Research (INSERM) and the French Blood Agency (EFS). The objective of observing children born within the same year is to analyze their physical and psychological health together with their living and environmental conditions. This large-scale study of children's development and socialization is the first of its kind in France. The collected data are now available to public and private research teams and many projects are underway in areas such as health, health environment and social sciences. In order to derive reliable confidence intervals for  finite population parameters such as totals or ratios, the ELFE sampling design has to be taken into account.

\noindent The ELFE sample is drawn according to a non-standard sampling design, called Cross-Classified Sampling (CCS), following Ohlsson (1996). It consists in drawing independently two samples from each component of a two-dimensional population. In the ELFE survey, a sample of maternity units and a sample of days are independently selected. CCS is peculiar in that the same sample of days is used for each of the selected maternity units, unlike two-stage sampling where independent sub-samples are selected inside the primary sampling units (see S\"arndal et al., 1992). Considering CCS as a particular two-phase sampling design is possible, and may prove to be useful as in section \ref{ssec:clt} below, but is fairly artificial. The two populations involved can be made of units with different natures (like days and maternities for the ELFE survey), leading to samples that are not part of one another but play a symmetric role. This sampling design appears in other contexts than the ELFE survey. Some examples include consumer price index surveys, as detailed in Dal\'en \& Olhsson (1995) for the Swedish survey, where outlets and items are sampled, and business surveys (Skinner, 2015), where businesses and products are sampled. Due to its particular properties, CCS deserves a specific attention. However, as noted by Skinner (2015), "the literature on the theory of cross-classified sampling is very limited". In particular, no general theory is derived under the finite population framework. While the papers by Vos (1964) and Ohlsson (1996) focus on simple random sampling without replacement, Skinner (2015) give some results under stratified without replacement simple random sampling and under with replacement unequal probability sampling.

\noindent In the present paper, we develop a general theory for estimation and variance estimation under CCS. The asymptotic normality of the Horvitz-Thompson estimator is derived under some mild conditions. A comparison with a two-stage sampling design is carried out in a general framework. We also raise an issue, not reported before, of possible negative values for Horvitz-Thompson and Yates-Grundy variance estimates. This problem occurs even in the simplest case of simple random sampling without replacement. Non-negative simplified variance estimators are therefore introduced. Conditions for their approximate unbiasedness are given under a design-based and a model-based approach. The properties of our variance estimators are evaluated through a small but realistic simulation study when estimating totals and ratios. Finally, an application to the ELFE data is detailed.


\section{\baselineskip=20pt Cross-classified sampling design} \label{sec:Varest}
\baselineskip=24pt


\subsection{\baselineskip=20pt Notations and Horvitz-Thompson estimation} \label{ssec:notation}

\noindent Keeping in mind the ELFE survey, we consider a population $U_M$ of $N_M$ maternities and a population $U_D$ of $N_D$ days. However, the developments below are completely general and may be applied to any populations $U_M$ and $U_D$. We will use the indexes $i$ and $j$ for the maternities, and the indexes $k$ and $l$ for the days. We consider a sampling design $p_M(\cdot)$ on the population $U_M$, leading to a sample $S_M$ of (average) size $n_M$, and a sampling design $p_D(\cdot)$ on the population $U_D$ leading to a sample $S_D$ of (average) size $n_D$. We assume that the two samples are selected independently. The cross-classified sampling design $p(\cdot)$ on the product population $U=U_M \times U_D$ is therefore defined as
    \begin{eqnarray*}
      p(s)=p_M(s_M) \times p_D(s_D) & \textrm{ for any } & s=s_M \times s_D \subset U_M \times U_D.
    \end{eqnarray*}

\noindent Let $\pi_i^M$ denote the probability that $i$ is selected in $S_M$, $\pi_{ij}^M$ denote the probability that units $i$ and $j$ are selected jointly in $S_M$, and let $\Delta_{ij}^M=\pi_{ij}^M-\pi_i^M \pi_j^M$. The quantities $\pi_k^D$, $\pi_{kl}^D$ and $\Delta_{kl}^D$ are similarly defined. We assume that the first and second-order inclusion probabilities are non-negative in each population. The probability for the couple $(i,k)$ to be selected in the product sample $S_M \times S_D$ is $\pi_i^M \pi_k^D$, and the probability for the couples $(i,k)$ and $(j,l)$ to be selected jointly in the product sample $S_M \times S_D$ is $\pi_{ij}^M \pi_{kl}^D$.

\noindent We are interested in some non-negative variable of interest with value $Y_{ik}$ for the maternity $i$ and the day $k$. The total $t_Y = \sum_{i\in U_M} \sum_{k\in U_D} Y_{ik}$ is then unbiasedly estimated by the Horvitz-Thompson (HT) estimator
    \begin{eqnarray} \label{estim:ht:prod}
    \hat{t}_{Y} & = & \sum_{i\in S_M} \sum_{k\in S_D} \frac{Y_{ik}}{\pi^M_i \pi^D_k}  = \sum_{i\in S_M} \sum_{k\in S_D} \check{Y}_{ik}  \;\; \mbox{ where }\;\; \check{Y}_{ik}=\frac{Y_{ik}}{\pi_i^M\pi_k^D}.
    \end{eqnarray}
Making use of the independence between $S_M$ and $S_D$, the variance of the HT-estimator is
    \begin{eqnarray} \label{eq:var}
    V_{CCS}\left(\hat{t}_{Y}\right) & = & \sum_{i,j \in U_M} \sum_{k,l \in U_D} \, \Gamma_{ijkl} \,\check{Y}_{ik}\check{Y}_{jl}
    \end{eqnarray}
where $\Gamma_{ijkl}=\pi_{ij}^M \pi_{kl}^D - \pi_{i}^M \pi_{j}^M \pi_{k}^D \pi_{l}^D$. The Sen(1953)-Yates-Grundy(1953) form
    \begin{eqnarray} \label{eq:varYG}
    V_{CCS}\left(\hat{t}_{Y}\right) & = & -\frac{1}{2}\sum_{(i,k) \neq (j,l) \in U_M \times U_D} \Gamma_{ijkl} \left(\check{Y}_{ik} - \check{Y}_{jl}\right)^2
    \end{eqnarray}
can  be used alternatively when both sampling designs are of fixed size.

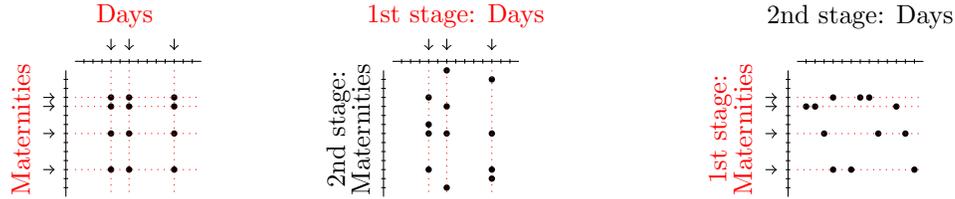
\begin{figure}[h] \footnotesize
\begin{minipage}{0.3\textwidth}
\begin{tikzpicture}[scale=0.6]
\raisebox{-0.6cm}{
\hspace{1cm}
\draw (1.3,7) node[text=red]{Days} ;
\draw (-1,4.5) node[text=red,rotate=90]{Maternities} ;

\draw (0.2,6)--(3,6);
\foreach \x in {0.4,0.6,...,3} {
\draw (\x,6.05cm) -- (\x,5.95cm) ;
}
\draw (0,5.8)--(0,3);
\foreach \x in {5.6,5.4,...,3} {
\draw (0.05cm,\x) -- (-0.05cm,\x) ;
}
\fill (1,5.2) circle (2pt) ; \fill (1.4,5.2) circle (2pt) ; \fill (2.4,5.2) circle (2pt) ;
\fill (1,5) circle (2pt) ; \fill (1.4,5) circle (2pt) ; \fill (2.4,5) circle (2pt) ;
\fill (1,4.4) circle (2pt) ; \fill (1.4,4.4) circle (2pt) ; \fill (2.4,4.4) circle (2pt) ;
\fill (1,3.6) circle (2pt) ; \fill (1.4,3.6) circle (2pt) ; \fill (2.4,3.6) circle (2pt) ;
\draw [<-] (1,6.25) -- (1,6.5); \draw [<-] (1.4,6.25) -- (1.4,6.5); \draw [<-] (2.4,6.25) -- (2.4,6.5);
\draw [<-] (-0.25,5.2) -- (-0.5,5.2); \draw [<-] (-0.25,5) -- (-0.5,5); \draw [<-] (-0.25,4.4) -- (-0.5,4.4); \draw [<-] (-0.25,3.6) -- (-0.5,3.6);

\draw[dotted,draw=red] (0.2,5.2)--(3,5.2);
\draw[dotted,draw=red] (0.2,5)--(3,5);
\draw[dotted,draw=red] (0.2,4.4)--(3,4.4);
\draw[dotted,draw=red] (0.2,3.6)--(3,3.6);
\draw[dotted,draw=red] (1,5.8)--(1,3);
\draw[dotted,draw=red] (1.4,5.8)--(1.4,3);
\draw[dotted,draw=red] (2.4,5.8)--(2.4,3);
}
\end{tikzpicture}
\end{minipage}
\hfill
\begin{minipage}{0.3\textwidth}
\begin{tikzpicture}[scale=0.6]
\raisebox{-0.6cm}{
\hspace{0cm}
\draw (1.6,2) node[text=red]{1st stage: Days} ;
\draw (-0.5,-0.5) node[text=,rotate=90]{Maternities} ;
\draw (-1,-0.5) node[text=,rotate=90]{2nd stage:} ;

\draw (0.2,1)--(3,1);
\foreach \x in {0.4,0.6,...,3} {
\draw (\x,1.05cm) -- (\x,0.95cm) ;
}
\draw (0,0.8)--(0,-2);
\foreach \x in {0.6,0.4,...,-2} {
\draw (0.05cm,\x) -- (-0.05cm,\x) ;
}

\fill (1,0.2) circle (2pt) ; \fill (1.4,0.8) circle (2pt) ; \fill (2.4,0.6) circle (2pt) ;
\fill (1,-0.4) circle (2pt) ; \fill (1.4,0) circle (2pt) ; \fill (2.4,-1.6) circle (2pt) ;
\fill (1,-0.6) circle (2pt) ; \fill (1.4,-0.6) circle (2pt) ; \fill (2.4,-0.6) circle (2pt) ;
\fill (1,-1.4) circle (2pt) ; \fill (1.4,-1.8) circle (2pt) ; \fill (2.4,-1.4) circle (2pt) ;

\draw [<-] (1,1.25) -- (1,1.5); \draw [<-] (1.4,1.25) -- (1.4,1.5); \draw [<-] (2.4,1.25) -- (2.4,1.5);
\draw[dotted,draw=red] (1,0.8)--(1,-2);
\draw[dotted,draw=red] (1.4,0.8)--(1.4,-2);
\draw[dotted,draw=red] (2.4,0.8)--(2.4,-2);
 }
\end{tikzpicture}
\end{minipage}
\hfill
\begin{minipage}{0.3\textwidth}
\begin{tikzpicture}[scale=0.6]
\raisebox{-0.6cm}{
\hspace{-.1cm}
\draw (1.6,4) node[text=]{2nd stage: Days} ;
\draw (-1,1.5) node[text=red,rotate=90]{Maternities} ;
\draw (-1.5,1.5) node[text=red,rotate=90]{1st stage:} ;
\draw (0.2,3)--(3,3);
\foreach \x in {0.4,0.6,...,3} {
\draw (\x,3.05cm) -- (\x,2.95cm) ;}
\draw (0,2.8)--(0,0);
\foreach \x in {2.6,2.4,...,0} {
\draw (0.05cm,\x) -- (-0.05cm,\x) ;}
\fill (1,2.2) circle (2pt) ; \fill (1.6,2.2) circle (2pt) ; \fill (1.8,2.2) circle (2pt) ;
\fill (0.6,2) circle (2pt) ; \fill (0.4,2) circle (2pt) ; \fill (2.4,2) circle (2pt) ;
\fill (0.8,1.4) circle (2pt) ; \fill (2,1.4) circle (2pt) ; \fill (2.6,1.4) circle (2pt) ;
\fill (1,0.6) circle (2pt) ; \fill (1.4,0.6) circle (2pt) ; \fill (2.8,0.6) circle (2pt) ;
\draw [<-] (-0.25,2.2) -- (-0.5,2.2); \draw [<-] (-0.25,2) -- (-0.5,2); \draw [<-] (-0.25,1.4) -- (-0.5,1.4); \draw [<-] (-0.25,0.6) -- (-0.5,0.6);
\draw[dotted,draw=red] (0.2,2.2)--(3,2.2);
\draw[dotted,draw=red] (0.2,2)--(3,2);
\draw[dotted,draw=red] (0.2,1.4)--(3,1.4);
\draw[dotted,draw=red] (0.2,0.6)--(3,0.6);
}
\end{tikzpicture}
\end{minipage}
\vspace{0.5cm}
\caption{Cross-classified sampling (left panel), two-stage sampling $DM$ with primary units in $U_D$ (central panel), two-stage sampling $MD$ with primary units in $U_M$ (right panel)} \label{fig:comp}
\end{figure}

\noindent Our set-up can be linked to the usual two-stage framework, by considering $U_M$ as a population of Primary Sampling Units (PSUs) and $U_D$ as a population of Secondary Sampling Units (SSUs), each maternity $i$ being associated to the same population of days. In case of two-stage sampling, denoted by $MD$, a first-stage sample $S_M$ is selected in $U_M$, and some second-stage samples $S_{i}$ are selected independently inside any $i \in S_M$. The variance of the HT-estimator is then
\begin{eqnarray}
V_{MD}\left(\hat{t}_{Y}\right) & = & V^{PSU}_{MD} \left(\hat{t}_Y\right) +V^{SSU}_{MD} \left(\hat{t}_Y\right) \label{var:ht:twost:1}
\end{eqnarray}
where
\begin{eqnarray}
V^{PSU}_{MD} \left(\hat{t}_Y \right) & = &  \sum_{i,j \in U_M} \sum_{k,l \in U_D} \Delta_{ij}^M \pi_k^D \pi_l^D \check{Y}_{ik} \check{Y}_{jl},\label{V:PSU:HT:1} \\
V^{SSU}_{MD} \left(\hat{t}_Y \right) & = &  \sum_{i \in U_M} \sum_{k,l \in U_D} \pi_i^M  \Delta_{kl}^D \check{Y}_{ik} \check{Y}_{il}. \label{V:SSU:HT:1}
\end{eqnarray}
Alternatively, we could consider $U_D$ as a population of PSUs and $U_M$ as a population of SSUs, each day $k$ being associated to the same population of maternities. In this case, the variance of the HT-estimator under two-stage sampling is
\begin{eqnarray}
V_{DM}\left(\hat{t}_{Y}\right) & = & V^{PSU}_{DM} \left(\hat{t}_Y\right) +V^{SSU}_{DM} \left(\hat{t}_Y\right) \label{var:ht:twost:2}
\end{eqnarray}
where
\begin{eqnarray}
V^{PSU}_{DM} \left(\hat{t}_Y \right) & = &  \sum_{k,l \in U_D} \sum_{i,j \in U_M} \Delta_{kl}^D \pi_i^M \pi_j^M \check{Y}_{ik} \check{Y}_{jl},\label{V:PSU:HT:2} \\
V^{SSU}_{DM} \left(\hat{t}_Y \right) & = &  \sum_{k \in U_D} \sum_{i,j \in U_M} \pi_k^D  \Delta_{ij}^M \check{Y}_{ik} \check{Y}_{il}. \label{V:SSU:HT:2}
\end{eqnarray}
The different features of CCS and two-stage sampling on a two-dimension population are illustrated on Figure \ref{fig:comp}.

\subsection{\baselineskip=20pt Variance decomposition for cross-classified sampling} \label{ssec:Var}

The covariance $\Gamma_{ijkl}$ may be written in several ways, leading to alternative variance decompositions. Plugging $\Gamma_{ijkl} = \pi_{kl}^D \Delta_{ij}^M + \pi_{ij}^M \Delta_{kl}^D - \Delta_{ij}^M \Delta_{kl}^D$ into (\ref{eq:var}) gives
\begin{eqnarray}
V_{CCS}\left(\hat{t}_{Y}\right) & = & V_{1} \left(\hat{t}_Y\right) +V_{2} \left(\hat{t}_Y\right) - V_{3} \left(\hat{t}_Y\right)\label{eq:var1-3}
\end{eqnarray}
where
\begin{eqnarray}
V_{1} \left(\hat{t}_Y\right) & = & \sum_{k,l \in U_D} \sum_{i,j \in U_M} \pi_{kl}^D \Delta_{ij}^M\, \check{Y}_{ik}\check{Y}_{jl}, \label{eq:var1MD}\\
V_{2} \left(\hat{t}_Y\right) & = & \sum_{i,j \in U_M} \sum_{k,l \in U_D} \pi_{ij}^M \Delta_{kl}^D\, \check{Y}_{ik}\check{Y}_{jl}, \label{eq:var1DM}\\
V_{3} \left(\hat{t}_Y\right) & = & \sum_{i,j \in U_M} \sum_{k,l \in U_D} \Delta_{ij}^M  \Delta_{kl}^D \check{Y}_{ik}\check{Y}_{jl}.\label{eq:var3}
\end{eqnarray}
Plugging $\Gamma_{ijkl} = \Delta_{ij}^M \pi_k^D \pi_l^D + \Delta_{kl}^D \pi_i^M \pi_j^M + \Delta_{ij}^M \Delta_{kl}^D$ into (\ref{eq:var}) gives
\begin{eqnarray}
V_{CCS}\left(\hat{t}_{Y}\right) & = & V^{PSU}_{MD} \left(\hat{t}_Y\right) +V^{PSU}_{DM} \left(\hat{t}_Y\right) +V_{3}^{} \left(\hat{t}_Y\right) \label{eq:var2+3}
\end{eqnarray}
and we have $V_{1} \left(\hat{t}_Y\right)=V^{PSU}_{MD} \left(\hat{t}_Y\right) +V_{3}^{} \left(\hat{t}_Y\right)$
and $V_{2} \left(\hat{t}_Y\right)=V^{PSU}_{DM} \left(\hat{t}_Y\right) +V_{3}^{} \left(\hat{t}_Y\right)$.
This second decomposition was originally derived by Ohlsson (1996). Other decompositions are possible, e.g. through an analysis of variance decomposition as for two-stage sampling.

\subsection{\baselineskip=20pt Comparison with two-stage sampling} \label{ssec:twostage}

\noindent From expressions (\ref{var:ht:twost:2}) and (\ref{eq:var2+3}), we obtain after some algebra that
    \begin{eqnarray} \label{comp:prod:twost:1}
      V_{CCS}\left(\hat{t}_{Y}\right)-V_{DM}\left(\hat{t}_{Y}\right)
      = \sum_{i,j \in U_M} \Delta_{ij}^M \sum_{k \neq l \in U_D} \pi_{kl}^D \check{Y}_{ik} \check{Y}_{jl}.
    \end{eqnarray}
In case of Poisson sampling (PO) inside $U_M$, the right-hand side in (\ref{comp:prod:twost:1}) is non-negative and CCS is thus less efficient than two-stage sampling. In case of fixed-size sampling inside $U_M$, equation (\ref{comp:prod:twost:1}) may be alternatively written as
    \begin{eqnarray} \label{comp:prod:twost:2}
      V_{CCS}\left(\hat{t}_{Y}\right)-V_{DM}\left(\hat{t}_{Y}\right)
      = \sum_{i \neq j \in U_M} \frac{(-\Delta_{ij}^M)}{2} \sum_{k \neq l \in U_D} \frac{\pi_{kl}^D}{\pi_{k}^D \pi_{l}^D} \left(\frac{Y_{ik}}{\pi_i^M}-\frac{Y_{jk}}{\pi_j^M}\right)\left(\frac{Y_{il}}{\pi_i^M}-\frac{Y_{jl}}{\pi_j^M}\right).
    \end{eqnarray}
If the so-called Sen-Yates-Grundy conditions are respected for $p_M$, the quantities $(-\Delta_{ij}^M)$ are non-negative. If $Y_{ik}$ is roughly proportional to the size of the maternity unit $i$, as can be expected for count variables, the quantities
$$\displaystyle \left(\frac{Y_{ik}}{\pi_i^M}-\frac{Y_{jk}}{\pi_j^M}\right)\left(\frac{Y_{il}}{\pi_i^M}-\frac{Y_{jl}}{\pi_j^M}\right)$$ will tend to be positive unless the inclusion probabilities $\pi_i^M$ are defined proportionally to some measure of size. CCS sampling would then be less efficient than two-stage sampling. This result is illustrated in section \ref{ssec:compsimu} on some simulated populations when both $p_M$ and $p_D$ are simple random sampling without replacement (SI) designs, and for different sample sizes.

\section{\baselineskip=20pt Variance estimation} \label{sec:HTvar}

\subsection{\baselineskip=20pt Design-unbiased variance estimation} \label{ssec:HTvar}

\noindent The HT variance estimator for $V_{CCS}\left(\hat{t}_{Y}\right)$ is
\begin{eqnarray} \label{eq:HT:var}
\hat{{V}}_{HT}\left(\hat{t}_{Y}\right) & = & \sum_{i,j \in S_M} \sum_{k,l \in S_D} \frac{\Gamma_{ijkl}}{\pi_{ij}^M \pi_{kl}^D}  \,\check{Y}_{ik}\check{Y}_{jl}.
\end{eqnarray}
It may be also derived from (\ref{eq:var1-3}), leading to the alternative writing
\begin{eqnarray}
\hat{{V}}_{HT}\left(\hat{t}_{Y}\right)
& = & \hat{V}_{1,HT}^{} \left(\hat{t}_Y\right) +\hat{V}_{2,HT}^{} \left(\hat{t}_Y\right) -\hat{V}_{3,HT}^{} \left(\hat{t}_Y\right) \label{eq:HT:var1-3}
\end{eqnarray}
where
\begin{eqnarray}
\hat{V}_{1,HT}^{} \left(\hat{t}_Y\right) &=& \sum_{i,j\in S_M} \sum_{k,l \in S_D} \frac{\Delta_{ij}^M}{\pi_{ij}^M}\, \check{Y}_{ik}\check{Y}_{jl}, \label{eq:estHT:var1MD} \\
\hat{V}_{2,HT}^{} \left(\hat{t}_Y\right) &=& \sum_{i,j\in S_M} \sum_{k,l \in S_D} \frac{\Delta_{kl}^D}{\pi_{kl}^D}\, \check{Y}_{ik}\check{Y}_{jl}, \label{eq:estHT:var1DM} \\
\hat{V}_{3,HT}^{} \left(\hat{t}_Y\right) &=& \sum_{i,j \in S_M} \sum_{k,l \in S_D} \frac{\Delta_{ij}^M}{\pi_{ij}^M}  \frac{\Delta_{kl}^D}{\pi_{kl}^D}  \,\check{Y}_{ik}\check{Y}_{jl}.\label{eq:estHT:var3}
\end{eqnarray}
If $p_M$ and $p_D$ are both Poisson sampling designs, this variance estimator is always non-negative. Otherwise, it may take negative values even if $p_M$ and $p_D$ are both
SI designs (denoted by SI$^2$) as illustrated in section \ref{ssec:simpsimu}. When $p_M$ and $p_D$ are both fixed-size sampling designs, we may alternatively consider the Yates-Grundy like variance estimator:
\begin{eqnarray}
\hat{V}_{YG}\left(\hat{t}_{Y}\right)
& = & \hat{V}_{1,YG}^{} \left(\hat{t}_Y\right) +\hat{V}_{2,YG}^{} \left(\hat{t}_Y\right) -\hat{V}_{3,YG}^{} \left(\hat{t}_Y\right)\label{eq:YG:var1-3}
\end{eqnarray}
where
\begin{eqnarray}
\hat{V}_{1,YG}^{} \left(\hat{t}_Y\right) &=&  -\frac{1}{2} \sum_{i \neq j \in S_M} \frac{\Delta_{ij}^M}{\pi_{ij}^M} \left(\frac{\hat{Y}_{i \bullet}}{\pi_i^M}-\frac{\hat{Y}_{j \bullet}}{\pi_j^M}\right)^2 \\
\hat{V}_{2,YG}^{} \left(\hat{t}_Y\right) &=&  -\frac{1}{2} \sum_{k \neq l \in S_D} \frac{\Delta_{kl}^D}{\pi_{kl}^D} \left(\frac{\hat{Y}_{\bullet k}}{\pi_k^D}-\frac{\hat{Y}_{\bullet l}}{\pi_l^D}\right)^2 \\
\hat{V}_{3,YG}^{} \left(\hat{t}_Y\right) &=& - \frac{1}{2}\sum_{(i,k) \neq (j,l) \in S_M \times S_D} \frac{\Delta_{ij}^M \Delta_{kl}^D}{\pi_{ij}^M \pi_{kl}^D} \left(\check{Y}_{ik} - \check{Y}_{jl}\right)^2
\end{eqnarray}
with $\hat{Y}_{\bullet k} = \sum_{i \in S_M} Y_{ik} / \pi_i^M$ is the estimated sub-total for the day $k$ and $\hat{Y}_{i \bullet} = \sum_{k \in S_D} Y_{ik} / \pi_k^D$ is the estimated sub-total for the maternity $i$.
It can be proved that $\hat{V}_{HT}\left(\hat{t}_{Y}\right)$ in (\ref{eq:HT:var1-3}) and $\hat{V}_{YG}\left(\hat{t}_{Y}\right)$ in (\ref{eq:YG:var1-3}) match term by term, when $p_M$ and $p_D$ are stratified simple random sampling designs. If both sampling designs satisfy the Sen-Yates-Grundy conditions (SYG), the terms $\hat{V}_{1,YG}^{} \left(\hat{t}_Y\right)$ and $\hat{V}_{2,YG}^{} \left(\hat{t}_Y\right)$ are non-negative. However, the term $\hat{V}_{3,YG}^{} \left(\hat{t}_Y\right)$ is usually non-negative, which may lead to negative values for $\hat{V}_{YG}^{} \left(\hat{t}_Y\right)$ as illustrated in the simulations of section \ref{ssec:simpsimu}. It is thus desirable to exhibit non-negative variance estimators with limited bias.

\subsection{\baselineskip=20pt Non-negative variance estimators} \label{ssec:Varestsimp}
\baselineskip=24pt

We consider the variance decomposition in (\ref{eq:var1-3}), and study the relative order of magnitude of the components. We make the following assumptions:
\begin{itemize}
  \item[H1:] There exist some constants $\alpha_1$ and $\alpha_2$ such that
    \begin{eqnarray*}
      \forall k \in U_D, ~~ \frac{1}{N_M} \sum_{i \in U_M} Y_{ik}^2 \leq \alpha_1, & \textrm{ and } & \forall i \in U_M, ~~ \frac{1}{N_D} \sum_{k \in U_D} Y_{ik}^2 \leq \alpha_2.
    \end{eqnarray*}
  \item[H2:] There exists some constants $\lambda_1>0$ and $\lambda_2>0$ such that
    \begin{eqnarray*}
      \forall k \in U_D, ~~ \pi_{k}^D \geq \lambda_1 \frac{n_D}{N_D}, & \textrm{ and } & \forall i \in U_M, ~~ \pi_{i}^M \geq \lambda_2 \frac{n_M}{N_M}.
    \end{eqnarray*}
  \item[H3:] There exist some constants $\gamma_1$ and $\gamma_2$ such that
    \begin{eqnarray*}
      \forall k \neq l \in U_D, ~~ \frac{N_D^2}{n_D} \sup_{k \neq l \in U_D} \left|\Delta_{kl}^D\right| \leq \gamma_1, & \textrm{ and } &
      \forall i \neq j \in U_M, ~~ \frac{N_M^2}{n_M} \sup_{i \neq j \in U_M} \left|\Delta_{ij}^M\right| \leq \gamma_2.
    \end{eqnarray*}
  \item[H4:] There exists some constant $\delta>0$ such that
    \begin{eqnarray*}
    V_{CCS}\left(\hat{t}_{Y}\right) & \geq & \delta N_M^2 N_D^2 \left( \frac{1}{n_M} + \frac{1}{n_D} \right).
    \end{eqnarray*}
\end{itemize}

\noindent It is assumed in (H1) that the variable $y$ has bounded moments of order 2 for each maternity $i$ and for each day $k$. Assumptions (H2) and (H3) are classical in survey sampling and are satistified for many sampling designs, see for example Cardot et al. (2013). It is assumed in (H4) that the variance of the HT-estimator under CCS sampling has the order $N_M^2 N_D^2 (n_M^{-1}+n_D^{-1})$. From assumptions (H1-H4), there exist some constants $C_1$, $C_2$ and $C_3$ such that
    \begin{eqnarray}
      \frac{V_{1}^{} \left(\hat{t}_Y\right)}{V_{CCS}\left(\hat{t}_{Y}\right)} & \leq & C_1 \,\frac{1}{1+n_M n_D^{-1}}, \label{order:V1:DM}\\
      \frac{V_{2}^{} \left(\hat{t}_Y\right)}{V_{CCS}\left(\hat{t}_{Y}\right)} & \leq & C_2 \,\frac{1}{1+n_D n_M^{-1}}, \label{order:V1:MD}\\
      \frac{V_{3} \left(\hat{t}_Y\right)}{V_{CCS}\left(\hat{t}_{Y}\right)} & \leq & C_3 \,\frac{1}{n_D n_M^{-1}+n_M n_D^{-1}} \label{order:V3}
    \end{eqnarray}
The proof is given in Appendix \ref{appen1}. It follows from (\ref{order:V1:DM})-(\ref{order:V3}) that if $n_D$ is large and $n_M$ is bounded, both $V_{2}^{} \left(\hat{t}_Y\right)$ and $V_{3} \left(\hat{t}_Y\right)$ are negligible and a non-negative simplified variance estimator can be derived by focusing on $V_{1}^{} \left(\hat{t}_Y\right)$ only. This leads to
    \begin{eqnarray} \label{vsimp:1}
      \hat{V}_{\text{SIMP1}} \left(\hat{t}_Y\right) & = & \hat{V}_{1,YG}^{} \left(\hat{t}_Y\right).
    \end{eqnarray}
If the sampling design $p_D$ satisfies the SYG conditions, this simplified estimator is always non-negative. In the particular $\mbox{SI}^2$ case, we obtain
    \begin{eqnarray} \label{vsimp:1:SI2}
      \hat{V}_{\text{SIMP1}} \left(\hat{t}_Y\right) & = & N_{M}^2 \left(\frac{1}{n_{M}}-\frac{1}{N_{M}}\right) s_{\hat{Y}_{\circ \bullet}}^2
    \end{eqnarray}
where
    \begin{eqnarray}
      s_{\hat{Y}_{\circ \bullet}}^2 &=& \frac{1}{n_{M}-1} \sum_{i \in S_{M}} \left( \hat{Y}_{i \bullet} - \frac{1}{n_{M}} \sum_{j \in S_{M}} \hat{Y}_{j \bullet}\right)^2 .
    \end{eqnarray}

\noindent Symmetrically, both $V_{1}^{} \left(\hat{t}_Y\right)$ and $V_{3} \left(\hat{t}_Y\right)$ may be seen as negligible if $n_M$ is large and $n_D$ is bounded. Another simplified variance estimator is thus
    \begin{eqnarray} \label{vsimp:2}
      \hat{V}_{\text{SIMP2}} \left(\hat{t}_Y\right) & = & \hat{V}_{2,YG}^{} \left(\hat{t}_Y\right).
    \end{eqnarray}
If the sampling design $p_M$ satisfies the SYG conditions, this estimator is non-negative. In the particular $\mbox{SI}^2$ case, we have
    \begin{eqnarray} \label{vsimp:2:SI2}
      \hat{V}_{\text{SIMP2}} \left(\hat{t}_Y\right) & = &  N_D^2 \left(\frac{1}{n_D}-\frac{1}{N_D}\right)s^2_{\hat{Y}_{\bullet\circ}}
    \end{eqnarray}
where
    \begin{eqnarray}
      s^2_{\hat{Y}_{\bullet\circ}} & = & \frac{1}{n_{D}-1} \sum_{k \in S_{D}} \left( \hat{Y}_{\bullet k} - \frac{1}{n_{D}} \sum_{l \in S_{D}} \hat{Y}_{\bullet l}\right)^2 .
    \end{eqnarray}

\noindent A third possible simplified variance estimator is
    \begin{eqnarray} \label{vsimp:3}
      \hat{V}_{\text{SIMP3}} \left(\hat{t}_Y\right) &=& \hat{V}_{\text{SIMP1}} + \hat{V}_{\text{SIMP2}} \nonumber \\
      & = & \hat{V}_{1,YG}^{} \left(\hat{t}_Y\right)+\hat{V}_{2,YG}^{} \left(\hat{t}_Y\right).
    \end{eqnarray}
This estimator is non-negative if both $p_D$ and $p_M$ satisfy the SYG conditions. It is approximately unbiased for $V_{CCS}\left(\hat{t}_{Y}\right)$ if $n_D$ is large and $n_M$ is bounded, or if $n_M$ is large and $n_D$ is bounded. In the particular $\mbox{SI}^2$ case
    \begin{eqnarray} \label{vsimp:3:SI2}
      \hat{V}_{\text{SIMP3}} \left(\hat{t}_Y\right) = N_{M}^2 \left(\frac{1}{n_{M}}-\frac{1}{N_{M}}\right) s_{\hat{Y}_{\circ \bullet}}^2 + N_D^2 \left(\frac{1}{n_D}-\frac{1}{N_D}\right)s^2_{\hat{Y}_{\bullet\circ}}.
    \end{eqnarray}
Similar formula can be easily derived in the case of stratified simple random sampling without replacement and will be used in Section 5.
\subsection{\baselineskip=20pt Relative bias under a superpopulation model} \label{sec:rb:mod}
\baselineskip=24pt

\noindent We consider the following superpopulation model
    \begin{eqnarray} \label{model}
    Y_{ik} =  \mu + \sigma_M U_i + \sigma_D V_k + \sigma_E W_{ik}
    \end{eqnarray}
where $U_i$, $V_k$ and $W_{ik}$ are independently generated according to a standard normal distribution. This is an analysis of variance model with two crossed random factors and without repetition. Let ``$E_m$" denote the expectation with respect to the model (\ref{model}) and ``$E_{p}$" denote the expectation with respect to the CCS design. For each simplified variance estimator
$\hat{{V}}_{\textrm{SIMP}i}$, $i=1,2,3$, the relative bias RB under the model and under the sampling design is defined by
    \begin{eqnarray}
    \textrm{RB}_{m,p} \left[\hat{{V}}_{\textrm{SIMP}i} \left(\hat{t}_Y\right) \right] = \frac{{E}_m \left\{{E}_p \left[\hat{{V}}_{\textrm{SIMPi}}\left(\hat{t}_Y\right)\right]  - V_{CCS} \left(\hat{t}_Y\right)\right\}}{{E}_m \left[ V_{CCS}\left(\hat{t}_Y\right) \right]}.
    \end{eqnarray}
In the $\mbox{SI}^2$ case, these relative biases are of the form
    \begin{eqnarray} \label{rb:mp:vsimp12}
    \textrm{RB}_{m,p} \left[\hat{{V}}_{\textrm{SIMP}i} \left(\hat{t}_Y\right) \right] &=&  - 1/ (1 + A_i)
    \end{eqnarray}
for $i=1$ and 2 and
    \begin{eqnarray} \label{rb:mp:vsimp3}
    \textrm{RB}_{m,p} \left[\hat{{V}}_{\textrm{SIMP}3} \left(\hat{t}_Y\right) \right] &=&   1/ (1 + A_3)
    \end{eqnarray}
for some positive constant $A_i$, $i=1,2,3,$ depending on $\sigma_M$, $\sigma_D$, $\sigma_E$ and $n_M$, $N_M$, $n_D$ and $N_D$, see equations (\ref{A1})-(\ref{A3}). Equations (\ref{rb:mp:vsimp12}) and (\ref{rb:mp:vsimp3}) imply that the two first simplified variance estimators are negatively biased while the third one is positively biased. Using the notations $r_M=\sigma^2_M/\sigma^2_E$, $r_D=\sigma^2_D/\sigma^2_E$, $f_M=n_M/N_M$ and $f_D=n_D/N_D$, we have
    \begin{eqnarray}
     A_1 & = & \frac{1 - f_M }{1 - f_D} \,\frac{n_D r_M +1 }{n_M r_D + f_M }, \label{A1}\\
     A_2 & = & \frac{1 - f_D }{1 - f_M}\, \frac{n_M r_D +1 }{n_D r_M + f_D }, \label{A2}\\
     A_3 & = & \frac{n_D r_M + f_D }{1 - f_D} + \frac{n_M r_D + f_M }{1 - f_M}. \label{A3}
    \end{eqnarray}
The bias of $\hat{{V}}_{\textrm{SIMP}1}$ increases from $-1$ to $0$ when $A_1$ increases, which occurs in particular when the ratio $r_M$ or the sample size $n_D$ increases. In other words, $\hat{{V}}_{\textrm{SIMP}1}$ will have a small bias under model (\ref{model}) if the variable of interest contains some maternity effect or if the number of sampled days is large enough. Symmetrically, $\hat{{V}}_{\textrm{SIMP}2}$ will have a small bias under model (\ref{model}) if the variable of interest contains some day effect or if the number of sampled maternities is large enough. The bias of $\hat{{V}}_{\textrm{SIMP}3}$ decreases from $1$ to $0$ when $A_3$ increases, which occurs in particular when $r_M$ or $r_D$ increases, or when $n_M$ or $n_D$ increases. In other words, $\hat{{V}}_{\textrm{SIMP}3}$ will have a small bias under model (\ref{model}) if the variable of interest contains some maternity or some day effect, or if the number of sampled days or the number of sampled maternities is large enough. The simulation study in section \ref{sec:Simu} supports these results, and confirm that the variance tends to be underestimated with $\hat{{V}}_{\textrm{SIMP}1}$ or $\hat{{V}}_{\textrm{SIMP}2}$, and overestimated with $\hat{{V}}_{\textrm{SIMP}3}$.

\subsection{A central-limit theorem} \label{ssec:clt}
\baselineskip=24pt

To produce confidence intervals with appropriate asymptotic coverage, it is of interest to state a central-limit theorem (CLT) for CCS. Roughly speaking, Theorem \ref{theo:clt} below states that if the HT-estimator follows a CLT under both sampling designs $p_D$ and $p_M$, then the HT-estimator also follows a CLT under CCS. It is derived almost directly from Theorem 2 in Chen and Rao (2007), and the proof is therefore omitted.

\begin{thm} \label{theo:clt}
Suppose that assumptions (H1)-(H4) hold. Suppose that
\begin{itemize}
  \item[H5] $\sigma_{1}^{-1} V_1 \rightarrow_{\mathcal{L}} \mathcal{N}(0,1)$, where $\rightarrow_{\mathcal{L}}$ stands for the convergence in distribution under the sampling-design, with
    \begin{eqnarray} \label{theo:clt:eq1}
      V_1 = \frac{1}{N} \left(\sum_{i \in S_M} \frac{Y_{i\bullet}}{\pi_i^M} - \sum_{i \in U_M} Y_{i\bullet} \right) & \textrm{ and } & \sigma_{1}^2=V(V_1).
    \end{eqnarray}
  \item[H6] $\sup_t |P(\sigma_{2}^{-1} U_1 \leq t|S_M)-\Phi(t)|=o_p(1)$, where $\Phi$ is the cumulative distribution function of the standard normal distribution, and where
    \begin{eqnarray} \label{theo:clt:eq2}
      U_1 = \frac{1}{N} \sum_{i \in S_M} \frac{1}{\pi_i^M} (\hat{Y}_{i\bullet}-Y_{i\bullet}) & \textrm{ and } & \sigma_{2}^2=V(U_1|S_M).
    \end{eqnarray}
  \item[H7] $\sigma_{1}^2/\sigma_{2}^2 \rightarrow_{P} \gamma^2$, where $\rightarrow_{P}$ stands for the convergence in probability under the sampling-design.
\end{itemize}
Then
    \begin{eqnarray} \label{theo:clt:eq3}
      \frac{N^{-1}(\hat{t}_Y-t_Y)}{\sqrt{\sigma_{1n}^2+\sigma_{2n}^2}} & \rightarrow_{\mathcal{L}} & \mathcal{N}(0,1).
    \end{eqnarray}
\end{thm}


\noindent For illustration, we consider the particular case when $p_D$ and $p_M$ are both SI designs. Suppose that (H2)-(H4) hold, and that (H1) is strengthened to
\begin{itemize}
  \item[H1b:] There exists $\delta>0$ and some constants $\alpha_1$ and $\alpha_2$ such that
    \begin{eqnarray*}
      \forall k \in U_D ~~ \frac{1}{N_M} \sum_{i \in U_M} Y_{ik}^{2+\delta} \leq \alpha_1, & \textrm{ and } & \forall i \in U_M ~~ \frac{1}{N_D} \sum_{k \in U_D} Y_{ik}^{2+\delta} \leq \alpha_2.
    \end{eqnarray*}
\end{itemize}
Then by using the CLT in Hajek~(1961), the assumption (H5) can be shown to hold. By mimicking the proof of Lemma 2 in Chen and Rao (1997), the assumption (H6) can be shown to hold as well.

\section{\baselineskip=20pt Simulations } \label{sec:Simu}
\baselineskip=24pt

\noindent In this Section, two artificial populations are first generated using the superpopulation model (\ref{model}). In Section \ref{ssec:compsimu}, CCS is compared with two stage sampling in terms of variance, which illustrates the results in Section \ref{ssec:twostage}. A Monte Carlo experiment is then presented in Section \ref{ssec:simpsimu}, and the variance estimators introduced in Section \ref{sec:HTvar} are compared for the estimation of a total. Some attention is paid to the issue of negative values for the unbiased variance estimator. In Section \ref{ssec:simu:ratio}, two other populations with two variables of interest for each are generated. We focus on variance estimation for a ratio, making use of the variance estimators introduced in Section \ref{sec:HTvar} with estimated linearized variables instead of the variable of interest. The results from Tables \ref{tab:total} and \ref{tab:ratio} are readily reproducible using the R code provided in the supplementary materials of the present paper.

\subsection{\baselineskip=20pt Comparison with two-stage sampling} \label{ssec:compsimu}
\baselineskip=24pt

\begin{figure}[h!] \scriptsize
\centering
\begin{tabular}{|cc|}
  \hline
  Population 1 & Population 2  \\
  $\sigma_M=5$, $\sigma_D=5$ & $\sigma_M=0.5$, $\sigma_D=5$  \\
  \includegraphics[scale=0.4]{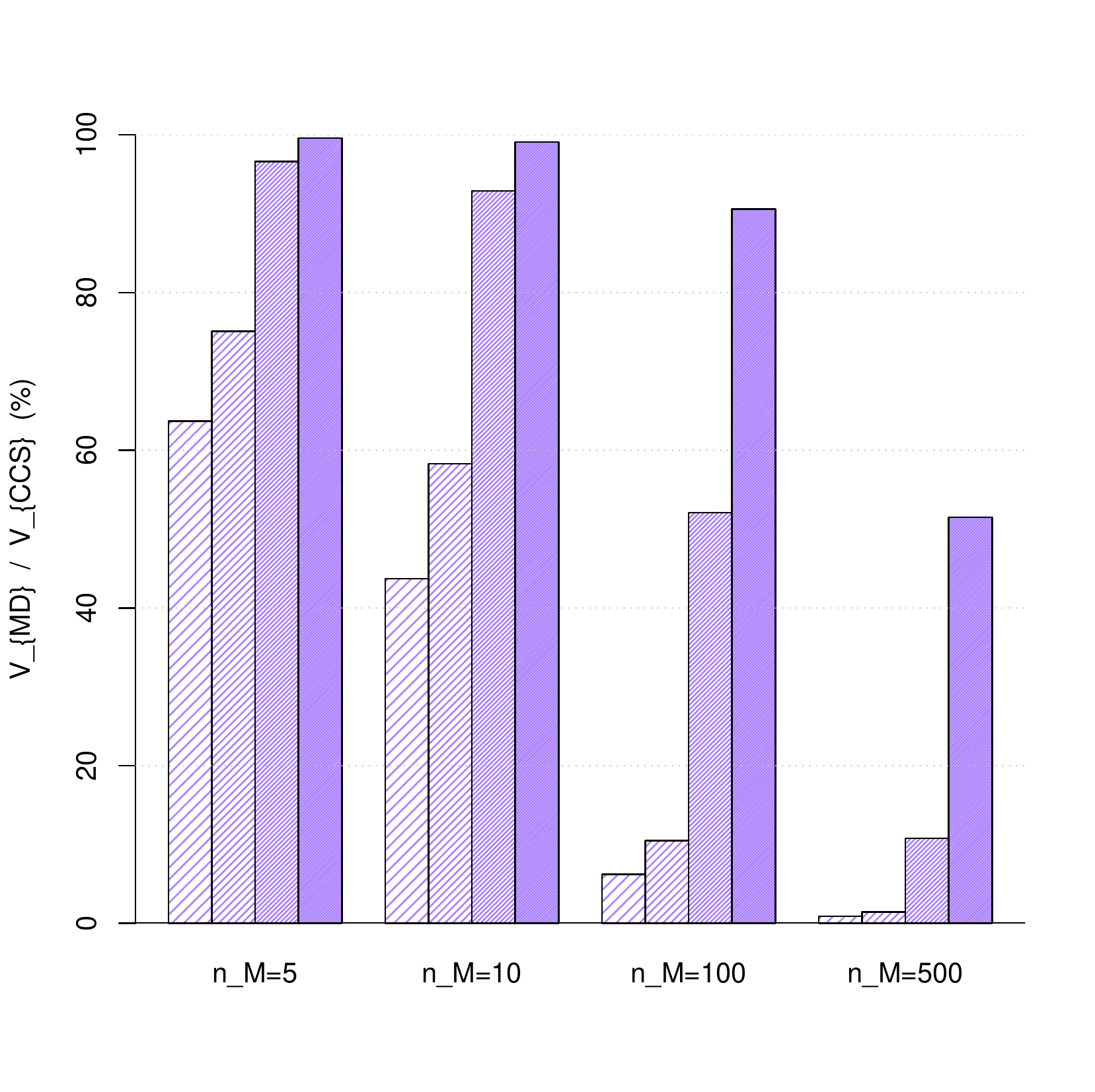}
& \includegraphics[scale=0.4]{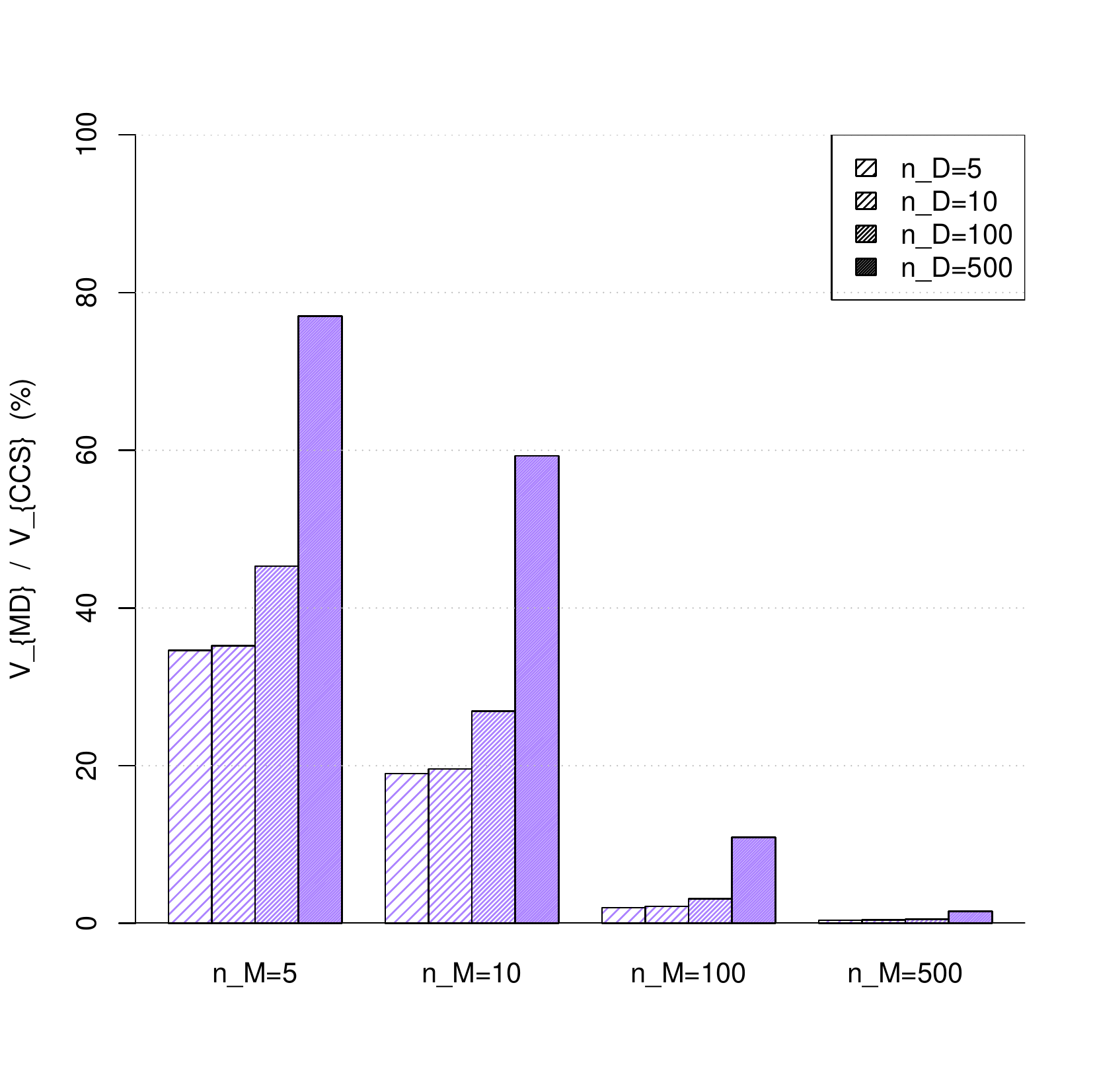}
 \\
  \hline
\end{tabular}
\caption{$V_{MD} \slash V_{CCS}$ ( \% ) for population 1 (left panel) and population 2 (right panel)} \label{graphCOMP}
\end{figure}

\noindent Two populations are generated according to model (\ref{model}), with $N_M=1000$ maternities and $N_D=1000$ days for each population, and with $\mu=200$ and $\sigma_E=5$. Equal random effects $\sigma_M=\sigma_D=5$ are used for population 1, while we use $\sigma_M=0.5$ and $\sigma_D=5$ for population 2. For each population, the $\mbox{SI}^2$ sampling design is used, with sample sizes equal to $5$, $10$, $100$ and $500$. The ratios $V_{MD} \slash V_{CCS}$ between the variance under two-stage sampling and the variance under CCS are computed, and plotted in percentage on Figure \ref{graphCOMP}. A ratio smaller than $100$ indicates that two-stage sampling is more accurate than CCS, which holds true in all cases considered in our experiment.

\noindent The ratio increases with $n_D$ and decreases when $n_M$ increases. Also, it can be observed that the ratio decreases with $\sigma_M$. This impact of the maternity effect is noticeable, and illustrates the substantial loss in accuracy induced by using a CCS instead of a two-stage sampling design if the maternity effect is small. Similar conclusions could be derived when computing the ratio $V_{DM} \slash V_{CCS}$.

\subsection{\baselineskip=20pt Variance estimation for a total} \label{ssec:simpsimu}
\baselineskip=24pt

\noindent We consider the two artificial populations generated as described in Section \ref{ssec:compsimu}. For each population, the $\mbox{SI}^2$ sampling design is used, with sample sizes equal to $5$, $10$, $100$ and $500$, and the sample selection is repeated $B=10,000$ times. For each sample $b=1,\ldots,B$, we compute the estimate $\hat{t}_Y^{(b)}$ of the total $t_Y$. The unbiased variance estimator $\hat{V}_{}^{(b)}$ and
the simplified variance estimators $\hat{V}_{\text{SIMP1}}^{(b)}$, $\hat{V}_{\text{SIMP2}}^{(b)}$, $\hat{V}_{\text{SIMP3}}^{(b)}$ are also computed for
$\hat{t}_Y^{(b)}$.

\noindent For each variance estimator $\hat{V}$, we compute the Monte Carlo Percent Relative Bias
 \begin{eqnarray*}
 \mbox{RB}_{\mbox{\sc mc}}(\hat{V}) = 100 \times \frac{B^{-1} \sum_{b=1}^B \hat{V}^{(b)}-V}{V},
 \end{eqnarray*}
where the true variance $V$ was approximated through an independent set of $50,000$ simulations. The number (\#NEG) of negative variance estimators $\hat{V}_{}^{(b)}$ is also computed.

\noindent The results are reported in Table \ref{tab:total}. The variance estimator $\hat{{V}}$ is almost unbiased in all situations, as expected. However,
this variance estimator is prone to negative values with small sample sizes when the value of $\sigma_M$ and/or the value of $\sigma_D$ is small as compared to
$\sigma_E$. The problem vanishes when the sample sizes increase. We now turn to the simplified variance estimators. The relative bias of
$\hat{{V}}_{\text{SIMP1}}$ decreases when $n_D$ increases or when $n_M$ decreases, and when $\sigma_M$ increases or when $\sigma_D$ decreases. This supports
the findings in Section \ref{sec:rb:mod}. Symmetrical conclusions are drawn for the relative bias of $\hat{{V}}_{\text{SIMP2}}$. Turning to
$\hat{{V}}_{\text{SIMP3}}$, we note that the relative bias decreases when either $\sigma_M$ or $\sigma_D$ increases. This variance estimator is therefore
advisable in all cases but those where there is no maternity nor day effect.

\begin{table}[htbp]
  \centering \scriptsize
    \begin{tabular}{l|rrrrr|rrrrr}
    \hline
    $n_M$ & 5     & 10    & 10    & 100   & 500   & 5     & 10    & 10    & 100   & 500    \\
    \hline
    $n_D$ & 5     & 10    & 100   & 100   & 500   & 5     & 10    & 100   & 100   & 500    \\
    \hline
    \hline
    ${\sigma}_M$ & \multicolumn{5}{c|}{5}  & \multicolumn{5}{c}{50} \\
    ${\sigma}_D$ & \multicolumn{5}{c|}{5}  & \multicolumn{5}{c}{5} \\
    \hline
    $\mbox{RB}_{\mbox{\sc mc}}\left(\hat{V}\right)$ & 1 & -1 & 2 & 0 & -0     & 1 & -1 & 1 & 0 & 0 \\
    \#NEG & 6 & 0 & 0 & 0 & 0  & 0  & 0 & 0 & 0 & 0 \\
    $\mbox{RB}_{\mbox{\sc mc}}\left(\hat{V}_{SIMP1}\right)$ & -43 & -47 & -6    & -49 & -49  & - & -2 & 1 & -1 & -1 \\
    $\mbox{RB}_{\mbox{\sc mc}}\left(\hat{V}_{SIMP2}\right)$ & -46 & -50 & -91 & -51 & -51  & -99  & -99  & -100   & -99   & -99 \\
    $\mbox{RB}_{\mbox{\sc mc}}\left(\hat{V}_{SIMP3}\right)$ & 11 & 3 & 2 & 1 & -0 & 1 & -0 & 1 & 0 & 0 \\
    \hline
    \hline
    ${\sigma}_M$ & \multicolumn{5}{c|}{0.5}   & \multicolumn{5}{c}{0.5} \\
    ${\sigma}_D$ & \multicolumn{5}{c|}{5}     & \multicolumn{5}{c}{0.5} \\
    \hline
    $\mbox{RB}_{\mbox{\sc mc}}\left(\hat{V}\right)$  & 1 & -1 & 0 & 1 & -1   & 1 & -1 & 2 & -0 & -0 \\
    \#NEG  & 91  & 0  & 0  & 0  & 0  & 1393 & 298 & 0  & 0  & 0 \\
    $\mbox{RB}_{\mbox{\sc mc}}\left(\hat{V}_{SIMP1}\right)$ & -82 & -90 & -81 & -98 & -99 & -4 & -9 & -3 & -34 & -47 \\
    $\mbox{RB}_{\mbox{\sc mc}}\left(\hat{V}_{SIMP2}\right)$ & -1 & -2 & -10 & -0 & -2 & -5 & -10 & -52 & -36 & -49 \\
    $\mbox{RB}_{\mbox{\sc mc}}\left(\hat{V}_{SIMP3}\right)$ & 18 & 8 & 9 & 2 & -0 & 90 & 81 & 45 & 29 & 4 \\
    \hline
    \end{tabular}
    \caption{Comparison between variance estimators for a  total}
  \label{tab:total}
\end{table}

\subsection{\baselineskip=20pt Variance estimation for a ratio} \label{ssec:simu:ratio}
\baselineskip=24pt

\noindent We now consider variance estimation for a ratio. Two populations are generated with $N_M=1000$ maternities and $N_D=1000$ days. In each population, two count variables are generated so as to mimic the data encountered in the ELFE survey. More precisely, we first generate an auxiliary variable $Z_{ik}$
according to model (\ref{model}) with $\mu=200$, $\sigma_E=\sigma_D=5$, and $\sigma_M=5$ or $50$. The first variable of interest $X_{ik}$ is generated
according to a Poisson distribution with parameter $Z_{ik}$. The second variable of interest $Y_{ik}$ is generated according to a binomial distribution with
parameters $X_{ik}$ and $p_{ik}$. We consider two cases: (i) equal probabilities with $p_{ik}=0.3$; (ii) unequal probabilities with $\logit(p_{ik})=\beta Z_{ik}$, where $\beta$ was chosen so that the average probability is approximately 0.3. Note that $Y_{ik}$ follows a Poisson distribution with parameter $p_{ik}Z_{ik}$.

\noindent The reason for this generating process is that some variable of interest $X_{ik}$, like the number of births in the ELFE survey, may contain some
maternity and/or day effect which is reflected in the way $Z_{ik}$ is generated. On the other hand, some maternity and/or day effect may also be contained in some other variable of interest $Y_{ik}$, like the number of births per caesarean. Such effects may be either similar to those for $X_{ik}$ like with pattern (i), or may occur differently like with pattern (ii).

\noindent For each population, the $\mbox{SI}^2$ sampling design is used, with sample sizes equal to $5$, $10$, $100$ and $500$, and the sample selection is
repeated $B=10,000$ times. For each sample $b=1,\ldots,B$, we compute the substitution estimator $\hat{R}^{(b)}={\hat{t}_Y^{(b)}} \slash {\hat{t}_X^{(b)}}$
of the ratio $R={t_Y} \slash {t_X}$. The variance estimator $\hat{V}_{}^{(b)}$ and the simplified variance estimators $\hat{V}_{\text{SIMP1}}^{(b)}$,
$\hat{V}_{\text{SIMP2}}^{(b)}$, $\hat{V}_{\text{SIMP3}}^{(b)}$ are also computed for $\hat{t}_Y^{(b)}$, where the variable of interest $Y_{ik}$ is replaced with
the estimated linearized variable of the ratio.

\noindent The results are reported in Table \ref{tab:ratio}. The variance estimator $\hat{{V}}$ is almost unbiased in all situations, as expected,
but is prone to negative values even when the maternity or day effect is small. We now turn to the relative bias for the simplified variance estimators. With
pattern (i), the situation is much different to that when a total is estimated, since the relative bias of $\hat{{V}}_{\text{SIMP3}}$ is much
larger than for the other two simplified estimators. This can be explained as follows: when the probabilities $p_{ik}$ are uniform, both $Y_{ik}$ and $X_{ik}$
contain the same maternity and day effect, but these effects wear off in the linearized variable. Whatever the values of $\sigma_M$ and $\sigma_D$ are, the situation is therefore comparable to that observed in the bottom right cell of Table \ref{tab:total}. With pattern (ii), the probabilities $p_{ik}$ depend on
$i$ and $k$, leading potentially to some remaining maternity and/or day effect in the linearized variable. In such situation, which seems more realistic in
practice, the relative bias of $\hat{{V}}_{\text{SIMP1}}$ and $\hat{{V}}_{\text{SIMP2}}$ increase when $\sigma_M$ or $\sigma_D$ increase, while the
relative bias of $\hat{{V}}_{\text{SIMP3}}$ decreases.

\begin{table}[htbp]
  \centering \scriptsize
    \begin{tabular}{l|l|rrrrr|rrrrr}
    \hline
    &$n_M$ & 5     & 10    & 10    & 100   & 500   & 5     & 10    & 10    & 100   & 500    \\
    \hline
    &$n_D$ & 5     & 10    & 100   & 100   & 500   & 5     & 10    & 100   & 100   & 500    \\
    \hline
    \hline
    &${\sigma}_M$ & \multicolumn{5}{c|}{5}  & \multicolumn{5}{c}{50} \\
    &${\sigma}_D$ & \multicolumn{5}{c|}{5}  & \multicolumn{5}{c}{5} \\
    \hline
    \hline
    Case (i) & $\mbox{RB}_{\mbox{\sc mc}}\left(\hat{V}\right)$ & -0 & -1 & -1 & 0 & -0 & -2 & -1 & -1 & 0 & 1 \\
    \multirow{4}{*}{$p_{ik}=0.3$}
    &\#NEG & 1645 & 484 & 14 & 0     & 0     & 1656 & 499 & 12 & 0 & 0 \\
    &$\mbox{RB}_{\mbox{\sc mc}}\left(\hat{V}_{SIMP1}\right)$ & -1 & -1 & -2 & -10 & -37 & -1 & -1 & -1 & -8 & -32 \\
    &$\mbox{RB}_{\mbox{\sc mc}}\left(\hat{V}_{SIMP2}\right)$ & 0 & -2 & -10 & -8 & -30 & -2 & -1 & -9 & -8 & -31 \\
    &$\mbox{RB}_{\mbox{\sc mc}}\left(\hat{V}_{SIMP3}\right)$ & 99 & 96 & 89 & 82 & 33 & 97 & 98 & 90 & 84 & 37 \\
    \hline
    \hline
    Case (ii) & $\mbox{RB}_{\mbox{\sc mc}}\left(\hat{V}\right)$ & 0 & -1 & 2 & 0 & -0 & -4 & -3 & -1 & -0 & 0 \\
    \multirow{4}{*}{$p_{ik} = \frac{e^{\beta Z_{ik}}}{1 + e^{\beta Z_{ik}}}$}
    &\#NEG & 1351  & 235 & 0     & 0     & 0     & 67 & 0     & 0     & 0     & 0  \\
    &$\mbox{RB}_{\mbox{\sc mc}}\left(\hat{V}_{SIMP1}\right)$ & -7 & -13 & -4 & -39 & -48 & -5 & -4  & -1 & -1 & -1  \\
    &$\mbox{RB}_{\mbox{\sc mc}}\left(\hat{V}_{SIMP2}\right)$ & -6 & -14 & -61 & -40 & -49 & -87 & -93 & -99 & -98 & -99 \\
    &$\mbox{RB}_{\mbox{\sc mc}}\left(\hat{V}_{SIMP3}\right)$ & 87 & 73 & 35 & 22 & 3 & 8 & 3 & -0 & 0 & 0 \\
    \hline
    \end{tabular}%
    \caption{Comparison between variance estimators for a ratio}
  \label{tab:ratio}%
\end{table}%

\section{\baselineskip=20pt Application to the ELFE survey} \label{sec:ELFE}
\baselineskip=24pt

\noindent ELFE is the first longitudinal study of its kind in France, tracking children from birth to adulthood. This cohort comprises more than 18,000 children
whose parents consented to their inclusion. The population of inference consists of babies born during 2011 in French maternities, excluding very premature
infants. It is a two-dimensional population with 544 maternities as spatial units and 365 days as time units. The crossing of one day and one maternity represents a cluster of infants.

\noindent The sample is obtained by CCS, where days and maternities are selected independently with selected families surveyed shortly after birth in $320$ metropolitan maternities and during $25$ days for one year. The sample selection for maternities may be modeled as stratified simple random sampling (STSI), the
population of maternities being divided into five strata of equal size. The allocation per stratum is proportional to the number of deliveries recorded in 2008.
The sample selection for days may be modeled as STSI, with four strata associated to the four seasons of 2011. The sample sizes inside strata are provided in
Tables \ref{tab:stsiMbis} and \ref{tab:stsiDbis}. \\

\begin {table}[h!]
\begin{center} \footnotesize
   \begin{tabular}{|l| c| c|}
    \hline
    Strata  & Strata size  &  Sample size      \\
    $g$ &  $N_{Mg}$ &   $n_{Mg}$     \\
    \hline
    1 & 108 & 21  \\     2 & 108 & 41  \\    3 & 109 & 55  \\     4 & 108 & 80  \\     5 & 111 & 90  \\
    \hline
    Total & 544 & 287\\
    \hline
    \end{tabular}
 \caption {Population and sample strata sizes for the maternities design $p_M$.} \label{tab:stsiMbis}
\end{center}
\end {table}
\begin {table}[h!]
\begin{center} \footnotesize
   \begin{tabular}{|l| c| c|}
    \hline
    Strata  & Strata size  &  Sample size      \\
    $h$ &  $N_{Dh}$ &   $n_{Dh}$     \\
    \hline
    1 & 91 & 4  \\     2 & 91 & 6  \\    3 & 91 & 7  \\     4 & 92 & 8  \\
    \hline
    Total & 365 & 25 \\
    \hline
    \end{tabular}
     \caption {Population and sample strata sizes for the days design $p_D$.} \label{tab:stsiDbis}
\end{center}
\end{table}

\noindent In this Section, we aim at illustrating the results previously obtained on a real data set. Some aspects of the ELFE survey, like the non-response
issue or the calibration step, deserve a specific attention but are beyond the scope of the present paper and are therefore not considered. In particular, the
ELFE survey is prone to several levels of non-response, since some sampled maternities and some families refused to participate either for some specific days
or for the whole period. In the present study, the sample of respondents is viewed as the original sample and in particular, we consider only the 287
maternities that participate during the 25 days of survey. The calibration step is not taken into account. The results below are meant to illustrate our
theoretical results, but are not intended for use in other contexts.

\noindent We consider seven count variables from the ELFE survey. Some of them depend on the characteristics of the maternities (e.g., the spatial location), like the variable indicating whether the mother is followed by a midwife. Others are related to the days of the survey, like the variable indicating whether the birth occurred by caesarean. For each variable, the estimated total $\hat{t}_{Y}$ from equation (\ref{estim:ht:prod}), the estimated variance
$\hat{\mathbf V} \left(\hat{t}_{Y}\right)$ from equation (\ref{eq:HT:var1-3}) and the three simplified estimators are given in the upper part of Table
\ref{tab:ELFE}. Similar indicators are given in the bottom part of Table \ref{tab:ELFE} for ratios, when the totals of the variables of interest are divided by
the total number of births.

\begin{table}[htbp] \tiny
  \centering
    \begin{tabular}{|l|rrrrrrrr|}
    \hline
    \textbf{} & \textbf{Birth} & \textbf{Born by} & \textbf{Twins} & \textbf{Born} & \textbf{To have} & \textbf{To have a} & \textbf{To have a }  & \textbf{To have an} \\
    \textbf{} & \textbf{} & \textbf{Caesarean} & \textbf{} & \textbf{within} & \textbf{a mother} & \textbf{mother aged}   & \textbf{primiparous} & \textbf{immigrant} \\
        \textbf{} & \textbf{} & \textbf{} & \textbf{} & \textbf{marriage} & \textbf{followed by} & \textbf{between 18} & \textbf{mother}& \textbf{mother}\\
         \textbf{} & \textbf{} & \textbf{} & \textbf{} &  & \textbf{a midwife} & \textbf{and 25 years}& & \\
    \hline
    $\hat{t}_Y$ & 362924 & 33873 & 10187 & 160283 & 42337 & 43238 & 162316 & 44169 \\
    $\hat{\mathbf{V}} \left(\hat{t}_Y\right)$ & 7,6E+07 & 1,5E+07 & 5,3E+05 & 2,0E+07 & 3,9E+06 & 2,6E+06 & 1,5E+07 &  3,6E+06 \\
    \hline
    $\mbox{RD}\left(\hat{\mathbf{V}}_{\text{SIMP1}} \right)$ & -63,7 \% & -95,5 \% & -63,5 \% & -64,6 \% & -13,2 \% & -49,7 \% & -46,5 \% & -58,2 \% \\
    $\mbox{RD}\left(\hat{\mathbf{V}}_{\text{SIMP2}} \right)$ & -31,1 \% & -1,9 \%  & -13,3 \% & -29,7 \% & -76,3 \% & -35,2 \% & -41,4 \%  & -33,4 \% \\
    $\mbox{RD}\left(\hat{\mathbf{V}}_{\text{SIMP3}} \right)$ & 5,2 \%   & 2,6 \%   & 23,2 \%  & 5,7 \%   & 10,5 \%  & 15,1 \%  & 12,2 \%  &  8,4 \%  \\
\hline \hline
    $\hat{R}$ & 1,00  & 0,09  & 0,03  & 0,44  & 0,12  & 0,12   &  0,45 &  0,12 \\
    $\hat{\mathbf{V}} \left(\hat{R}\right)$ &       & 7,9E-05 & 2,8E-06 & 2,4E-05 & 2,5E-05 & 1,2E-05  &  3,0E-05  & 1,6E-05   \\
    \hline
    $\mbox{RD}\left(\hat{\mathbf{V}}_{\text{SIMP1}} \right)$ &       & -96,2 \% & -51,0 \% & -31,0 \% & -7,9 \%  & -40,2 \%  & -69,3 \%  &  -49,2 \%  \\
    $\mbox{RD}\left(\hat{\mathbf{V}}_{\text{SIMP2}} \right)$ &       & -0,4 \%  & -17,0 \% & -44,7 \% & -80,5 \% & -35,5 \%  & -5,0 \%  &   -37,5 \% \\
    $\mbox{RD}\left(\hat{\mathbf{V}}_{\text{SIMP3}} \right)$ &       & 3,4 \%   & 31,9 \%  & 24,3 \%  & 11,5 \%  & 24,3 \%   &  25,7 \% &  13,3 \%  \\
    \hline
    \end{tabular}%
  \caption{Variance estimates of estimated total and ratio on some ELFE variables} \label{tab:ELFE}%
\end{table}%

\noindent The relative difference $RD$ between $\hat{V}_{\text{SIMP}}$ and the unbiased estimator $\hat{V}$ is
    \begin{eqnarray*}
    RD = \frac{\hat{V}_{\text{SIMP}} \left(\hat{t}_{Y\star}\right) - \hat{V} \left(\hat{t}_{Y\star}\right)}{\hat{V} \left(\hat{t}_{Y\star}\right)}.
    \end{eqnarray*}
Different behaviours may be observed for the variables of interest, depending on the maternity/day effect. For instance, the variable indicating whether the birth
occurred by caesarean contains an important day effect, and the RD of $\hat{{V}}_{\text{SIMP2}}$ is therefore small while that of
$\hat{{V}}_{\text{SIMP1}}$ is large. Symmetrically, the variable indicating whether the mother is followed by a midwife contains a small day effect as
compared to the maternity effect, and the RD of $\hat{{V}}_{\text{SIMP2}}$ is therefore large while that of
$\hat{{V}}_{\text{SIMP1}}$ is small. Also, we note that the RD of $\hat{{V}}_{\text{SIMP3}}$ is relatively stable for all variables when estimating a
total, which is an important feature in favour of this third simplified estimator. We note however that the absolute RD of $\hat{{V}}_{\text{SIMP3}}$ can
be large when estimating a ratio, which confirms the simulation results.

\section{Conclusion}

\noindent The present paper derives some general estimation theory for the cross-classified sampling design which was used in  the recent ELFE survey on childhood. The issue of possibly negative variance estimates arised even in case of simple random sampling without replacement. Alternative estimators to the usual Horvitz-Thompson and Yates-Grundy variance estimators are thus proposed, and proved to be non-negative under the usual Sen-Yates-Grundy conditions. The relative bias of the proposed variance estimators is derived for a superpopulation model. The behavior of these estimators is also investigated for totals and ratios on simulated data and on data extracted from the ELFE survey. Among the proposals, one variance estimator that leads to a slight overestimation of the variance in many cases, appears to be advisable. \\

\noindent Despite the present results and the recent paper by Skinner (2015), the cross-classified sampling design still deserves some attention. In particular, the treatment of non-response and the calibration problem should also be taken into account, and is currently under investigation.

\section{Supplementary Materials}

The three supplemental files are contained in a single archive (and can be obtained via a single download).

\begin{description}

\item[readme:] description of the supplemental files. (txt file)

\item[CodeR\_functions:] basic functions required to calculate estimators. (R file)

\item[CodeR\_Tables:] commands that calculate and display the results in Table 1 and Table 2 (call the CodeR\_functions). (R file)

\end{description}

\section{Bibliography}

\noindent Cardot, H., and Goga, C. and Lardin, P. (2013). \newblock{Uniform convergence and asymptotic confidence bands for model-assisted estimators of the
mean of sampled functional data.} \newblock{Electronic Journal of Statistics}, 7, 562-596. \\

\noindent Chen, J. and Rao, J.N.K. (2007). \newblock {Asymptotic normality under two-phase sampling designs.} \newblock {Statistica Sinica}, 17, 1047-1064. \\

\noindent Dal\'en, J. and Ohlsson, E. (1995). \newblock {Variance Estimation in the Swedish Consumer Price Index.} \newblock {Journal of Business \& Economic Statistics}, 13, No.3, 347-356. \\

\noindent Hajek, J. (1961). \newblock{Some extensions of the Wald-Wolfowitz-Noether theorem}. \newblock{The Annals of Mathematical Statistics}, 32, 506-523. \\

\noindent Ohlsson, E. (1996). \newblock {Cross-Classified Sampling.} \newblock {Journal of Official Statistics}, 12, No.3, 241-251. \\

\noindent S\"arndal, C.-E., Swensson, B. and Wretman, J.H. (1992). \newblock{Model Assisted Survey Sampling}. \newblock{New-York, Springer-Verlag.} \\

\noindent Sen, A.R. (1953). \newblock{On the estimate of the variance in sampling with varying probabilities.} \newblock{Journal of the Indian Society of Agricultural Statistics}, \textbf{5}, 119-127. \\

\noindent Skinner, C.J. (2015). \newblock{Cross-classified sampling: some estimation theory.} \newblock{Statistics and Probability Letters}, \textbf{104}, 163-168. \\

\noindent Vos, J. W. E. (1964). \newblock{Sampling in space and time.} \newblock{Review of the International Statistical Institute}, \textbf{32}, No. 3, 226-241. \\

\noindent Yates, F. and Grundy, P.M. (1953). \newblock{Selection without replacement from within strata with probability proportional to size.} \newblock{Journal of the Royal Statistical Society B}, \textbf{15}, 235-261.

\section{Appendix}
\subsection*{Proof of equations (\ref{order:V1:DM})-(\ref{order:V3})} \label{appen1}

We can rewrite
    \begin{eqnarray} \label{proof1:eq1}
      V_{1} \left(\hat{t}_Y\right) & = & \sum_{k \in U_D} \frac{V(\hat{Y}_{\bullet k})}{\pi_k^D} + \sum_{k \neq l \in U_D} \frac{\pi_{kl}^D}{\pi_{k}^D\pi_{l}^D} Cov(\hat{Y}_{\bullet k},\hat{Y}_{\bullet l}).
    \end{eqnarray}
We have
    \begin{eqnarray} \label{proof1:eq2}
      V(\hat{Y}_{\bullet k}) & = & \sum_{i \in U_M} (1-\pi_i^M) \frac{(Y_{ik})^2}{\pi_i^M}+ \sum_{i \neq j \in U_M} \frac{\pi_{ij}^M-\pi_i^M \pi_j^M}{\pi_i^M \pi_j^M} Y_{ik} Y_{jk}.
    \end{eqnarray}
From assumptions (H1), (H2) (H3) and Cauchy-Schwarz inequality, there exists some constant $C$ such that for any $k \in U_D$,
    \begin{eqnarray} \label{proof1:eq3}
      V(\hat{Y}_{\bullet k}) & \leq & C \frac{N_M^2}{n_M}.
    \end{eqnarray}
Also, from the Cauchy-Schwarz inequality, there exists some constant $C$ such that for any $k \neq l \in U_D$:
    \begin{eqnarray} \label{proof1:eq4}
      Cov(\hat{Y}_{\bullet k},\hat{Y}_{\bullet l}) & \leq & C \frac{N_M^2}{n_M}.
    \end{eqnarray}
From equation (\ref{proof1:eq3}) and assumption (H2), the first term in the right hand sum of (\ref{proof1:eq1}) is $O(N_D^2 N_M^2 n_M^{-1} n_D^{-1})$. From equation (\ref{proof1:eq4}) and assumptions (H2) and (H3), the absolute value of the second term in the RHS of (\ref{proof1:eq1}) is $O(N_D^2 N_M^2 n_M^{-1})$. Therefore, there exists some constant $C$ such that
    \begin{eqnarray} \label{proof1:eq5}
      V_{1} \left(\hat{t}_Y\right) & \leq & C \frac{N_D^2 N_M^2}{n_M}.
    \end{eqnarray}
We can prove similarly that there exists some constant $C$ such that
    \begin{eqnarray} \label{proof1:eq6}
      V_{2} \left(\hat{t}_Y\right) & \leq & C \frac{N_D^2 N_M^2}{n_D}.
    \end{eqnarray}
From equation (\ref{eq:var3}), the term $V_{3} \left(\hat{t}_Y\right)$ may be split into four terms according to the intersection of $\{i,j\}$ and $\{k,l\}$. From assumptions (H1)-(H3), it is easily shown that the absolute value of each of these four terms is $O(N_D^2 N_M^2 n_M^{-1} n_D^{-1})$. Therefore, there exists some constant $C$ such that
    \begin{eqnarray} \label{proof1:eq7}
      V_{3} \left(\hat{t}_Y\right) & \leq & C \frac{N_D^2 N_M^2}{n_M n_D}.
    \end{eqnarray}
Equations (\ref{order:V1:DM})-(\ref{order:V3}) follow immediately from equations (\ref{proof1:eq5})-(\ref{proof1:eq7}) and assumption (H4).

\end{document}